\documentclass[12pt,twoside]{article} 
\usepackage{amsmath,amssymb,latexsym,theorem,bbm,shapepar,natbib,epsfig}
\newfont{\msa}{msam10 scaled\magstep1}
\newfont{\ssmsa}{msam9}

\def\Sign{\mathop{\hbox{\rm sign}}}
    
\newcommand{\Var}{\mathrm{Var}}
\newcommand{\Cov}{\mathrm{Cov}}

\newcommand{\distr}{\stackrel{\scriptstyle{\cal D}}{\longrightarrow}}

\newcommand{\proofend}{\hfill$\square$}

\numberwithin{equation}{section}

\theorembodyfont{\em}
\newtheorem{Lem}{Lemma}[section]
\newtheorem{Thm}[Lem]{Theorem}

\newtheorem{Cor}[Lem]{Corollary}
\theorembodyfont{\rm}

\title{On the variances of a spatial unit root model}

\author{{\sc S\'andor Baran}\\
         Faculty of Informatics, University of Debrecen\\
         Egyetem t\'er 1, H--4032 Debrecen, Hungary}

\date{}

\begin{document}
\pagestyle{headings}

\maketitle

\begin{abstract}
Asymptotic properties of the variances of the
spatial autoregressive model \ $X_{k,\ell}
  =\alpha X_{k-1,\ell}+\beta X_{k,\ell-1}+\gamma X_{k-1,\ell-1}
   +\varepsilon_{k,\ell}$ \ are investigated in the unit root case,
   that is when the parameters are on the boundary of the domain of
   stability that forms a tetrahedron in \ $[-1,1]^3$. \ The limit
   of the variance of \ $n^{-\varrho}X_{[ns],[nt]}$ \ is determined,
   where on the interior of the faces of the domain of stability \ 
 $\varrho=1/4$,  \  on the  
edges \ $\varrho =1/2$, \ while on the vertices \ $\varrho =1$.  

\end{abstract}

\section{Introduction}
   \label{sec:sec1}
\markboth{\rm Baran}{\rm On the variances of a spatial
  unit root model} 

The analysis of spatial autoregressive models is of interest 
in many different fields of science such as 
 geography, geology, biology and agriculture. A detailed discussion of these
applications is given by \citet{br2} where the  
authors considered a special case of the so called 
unilateral autoregressive model having the form 
 \begin{equation}
    \label{brmod}
  X_{k_,\ell}
  =\sum_{i=0}^{p_1}\sum_{j=0}^{p_2}\alpha_{i,j}X_{k-i,\ell-j}
   +\varepsilon_{k,\ell}, 
  \qquad\alpha_{0,0}=0.
\end{equation}
A particular case of the  above model is the
doubly geometric spatial autoregressive process
 $$
  X_{k,\ell}
  =\alpha X_{k-1,\ell}+\beta X_{k,\ell-1}-\alpha\beta X_{k-1,\ell-1}
   +\varepsilon_{k,\ell},
 $$
 introduced by \citet{martin1}. This was the first spatial 
autoregressive model for which unstability has been studied.
It is, in fact, the simplest spatial model, since the product structure
 \ $\varphi(x,y)=xy-\alpha x-\beta y+\alpha\beta=(x-\alpha)(y-\beta)$ \ of
 its characteristic polynomial ensures that it can be considered as some kind
 of combination of two autoregressive processes on the line, and several
 properties can be derived by the analogy of one--dimensional autoregressive
 processes.
This model has been used by \citet{jain} in the study of image processing,
 by \citet{martin2}, \citet{cg}, \citet{br1} in agricultural trials and 
by \citet{tj1} in digital filtering. 

In the stable case when \ $|\alpha|<1$ \ and
 \ $|\beta|<1$, \ asymptotic normality of several estimators
 \ $(\widehat\alpha_{m,n},\widehat\beta_{m,n})$ \ of \
 $(\alpha,\beta)$ \ based on the observations
 \ $\{X_{k,\ell}:\text{$1\leq k\leq m$ \ and \ $1\leq\ell\leq n$}\}$ \ has been
 shown (e.g. \citet{tj2,tj3} or \citet{br3,br2}), namely,
 $$
  \sqrt{mn}
  \begin{pmatrix}
   \widehat\alpha_{m,n}-\alpha\\ 
   \widehat\beta_{m,n}-\beta
  \end{pmatrix}
  \distr{\mathcal N}(0,\Sigma_{\alpha,\beta})
 $$
as  \ $m,n\to\infty$ \ with \ $m/n\to\,\textup{constant}>0$ \ with some 
covariance matrix \ $\Sigma_{\alpha,\beta}$.
 
In the unstable  case when \ $\alpha=\beta=1$, \ in contrast
 to the classical first order autoregressive time series model, where 
the appropriately normed least squares estimator (LSE) of the autoregressive 
parameter converges to a fraction of functionals of the
standard Brownian motion (see e.g. \citet{phil} or \citet{cw}),
the sequence of Gauss--Newton estimators
 \ $(\widehat\alpha_{n,n},\widehat\beta_{n,n})$ \ of \
 $(\alpha,\beta)$ \ has been shown to 
 be asymptotically normal (see \citet{bhat1} and  \citet{bhat2}). 
In the unstable case \ $\alpha=1$, \ $|\beta|<1$ \ the LSE turns out
to be asymptotically normal again \citep{bhat1}.  

\citet{bpz1} discussed  a special case of the
model \eqref{brmod}, namely, when
 \ $p_1=p_2=1$, \ $\alpha_{0,1}=\alpha_{1,0}=:\alpha$ \ and \
 $\alpha_{1,1}=0$, which is the simplest spatial model, that can not
 be reduced somehow to autoregressive models on the line. This model
 is  stable  in case \ $|\alpha|<1/2$
 \ (see e.g. \citet{whittle}, \citet{besag} or \citet{br2}), and unstable  if
 \ $|\alpha|=1/2$. \ In \citet{bpz1} the asymptotic normality of the
LSE of the unknown parameter \ $\alpha$ \ 
is proved both in stable and unstable cases. The case  $p_1=p_2=1$, \
$\alpha_{1,0}=:\alpha, \ \alpha_{0,1}=:\beta$ \ and \
 $\alpha_{1,1}=0$ \ was studied by \citet{paul} and \citet{bpz2}. This model
 is  stable  in case \ $|\alpha|+|\beta|<1$
 \  and unstable  if \ $|\alpha|+|\beta|=1$ \ \citep{br2}.
 \citet{paul} determined the exact asymptotic behaviour of the
 variances of the process, while \citet{bpz2} proved the asymptotic
 normality of the LSE of the parameters both in
 stable and unstable cases.

In the present paper we study the asymptotic properties of a more 
complicated  special case of the model \eqref{brmod} with 
 \ $p_1=p_2=1$, \   $\alpha_{1,0}=:\alpha$, \ $\alpha_{0,1}=:\beta$ \ and \
$\alpha_{1,1}=:\gamma$. \ Our aim is to clarify the asymptotic
behaviour of the variances. The asymptotic results on the variances
(and covariances) 
help in finding the asymptotic properties of various estimators of the
autoregressive parameters (see e.g. \citet{bpz1,bpz2}). 

We consider the spatial autoregressive process
\ $\{X_{k,\ell}:k,\ell\in{\mathbb Z},\,k, \ell\geq0\}$ \ is defined as
 \begin{equation}\label{model}
  \begin{cases} X_{k,\ell}=
                     \alpha X_{k-1,\ell}+\beta X_{k,\ell-1}+\gamma
                     X_{k-1,\ell-1}+\varepsilon_{k,\ell}, 
                    & \text{for \ $k,\ell \geq 1$,}\\
                    X_{k,0}=X_{0,\ell}=0,  & \text{for \ $k,\ell\geq0$.}
                   \end{cases}
 \end{equation}
The model is stable if  \ $\ |\alpha |<1, \ |\beta |<1$ \ and \ $|\gamma|<1$, \
$|1+\alpha ^2-\beta ^2  -\gamma ^2|>2|\alpha +\beta\gamma|$ \ and \
$1-\beta ^2>|\alpha +\beta\gamma|$, \ and unstable on the boundary of
this domain \citep{br2} (see Figure \ref{fig:fig1}). Short calculation
shows that condition of stability means that \ $\
|\alpha |<1, \ |\beta |<1$ 
\ and \ $|\gamma|<1$, \ and inequalities 
\begin{equation*}
\alpha-\beta -\gamma <1, \qquad  -\alpha+\beta -\gamma <1, \qquad
-\alpha-\beta +\gamma<1, \qquad  \alpha+\beta +\gamma <1 
\end{equation*}
hold. 
Obviously, in case \ $\alpha\beta\gamma\geq 0$ \ the above set of conditions
reduces to \
$|\alpha|+|\beta|+|\gamma|<1$. \ If the model is not stable, 
one can distinguish three cases:
\begin{figure}[tb]
 \begin{center}
\leavevmode
\epsfig{file=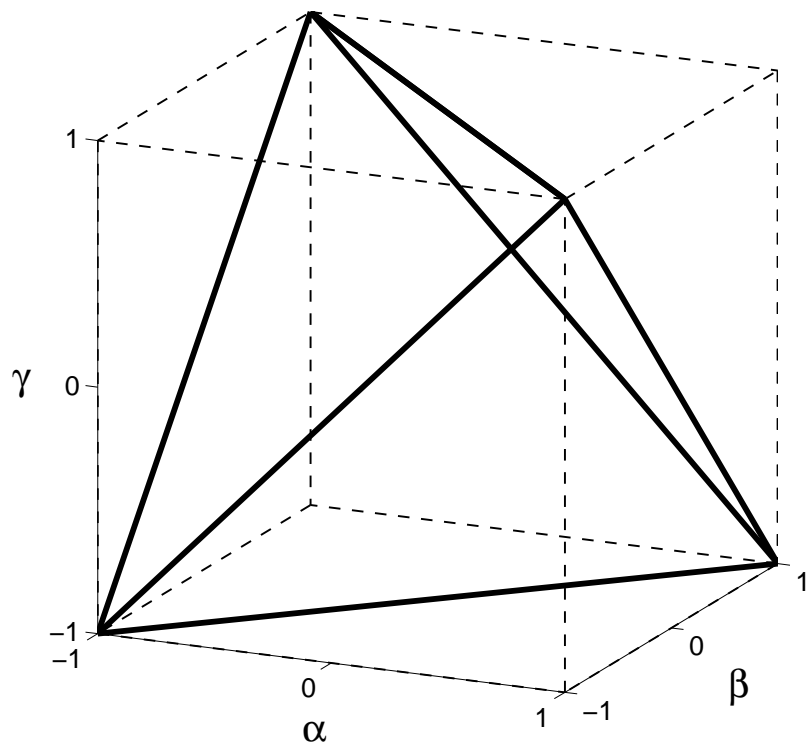,height=11cm}
\vskip1mm
\caption{The domain of stability of model \eqref{model}.}  
\label{fig:fig1}
  \end{center}  
\end{figure}

\noindent
{\sc Case A.} \ The parameters are in the interior of the faces of
the boundary of the domain of stability, i.e. \ $|\alpha|< 1, \
|\beta|< 1, \ |\gamma|< 1$ \ and one of the following equations is fulfilled
\begin{equation}
   \label{caseacond}
\alpha-\beta -\gamma =1; \qquad -\alpha+\beta -\gamma =1; \qquad -\alpha-\beta
+\gamma =1; \qquad \alpha+\beta +\gamma =1. 
\end{equation}
We remark that in case \ $\alpha\beta\gamma\geq 0$ \ the set of
equations \eqref{caseacond} is equivalent to \
$|\alpha|+|\beta|+|\gamma|=1$, 
\ while in case \ $\alpha\beta\gamma< 0$ \ to
\begin{equation*}
|\alpha|+|\beta| -|\gamma| =1 \quad \text{or} \quad
|\alpha|-|\beta|+|\gamma| =1  \quad \text{or} \quad -|\alpha|+|\beta|
+|\gamma| =1. 
\end{equation*}

\noindent
{\sc Case B.} \ The parameters are in the interior of the edges of
the boundary of the domain of stability, i.e. \ $\alpha\beta\gamma \leq 0$ \
and one of the following equations is fulfilled
\begin{equation*}
|\alpha|=1 \ \ \text{and} \ \ |\beta|=|\gamma|<1; \quad |\beta|=1 \ \
\text{and} \ \ |\alpha|=|\gamma|<1; \quad |\gamma|=1 \ \ \text{and} \ \ 
|\alpha|=|\beta|<1. 
\end{equation*}
Observe that in each of the above three cases exactly two of the defining
equations \eqref{caseacond} of Case A are satisfied. In this way Case B can be
considered as an extension of Case A to the situation when \
$\alpha\beta\gamma \leq 0$ \ and one of the
parameters equals \ $\pm 1$, \ while the other two parameters have
absolute values less than one.

\noindent
Further, observe that in the first two cases \ $\gamma=-\alpha \beta$, \ so we
obtain spacial cases of the doubly geometric model. If \ $|\alpha|=1, \
|\beta|=|\gamma| \leq 1$  \ then for \ $k\in{\mathbb N}$ \ the
difference \ 
$\Delta_{1,\alpha}X_{k,\ell}:=X_{k,\ell}-\alpha X_{k-1,\ell}$ \ is a
classical AR(1) process, i.e. \ $\Delta_{1,\alpha}X_{k,\ell}=\beta
\Delta_{1,\alpha}X_{k,\ell-1} +\varepsilon_{k,\ell}$. \ Similarly, if \
 $|\beta|=1, \ 
|\alpha|=|\gamma|\leq 1$  \ then \ $\Delta_{2,\beta}X_{k,\ell}=\alpha
\Delta_{2,\beta}X_{k-1,\ell} +\varepsilon_{k,\ell}$, \ where \
$\Delta_{2,\beta}X_{k,\ell}:=X_{k,\ell}-\beta X_{k,\ell-1}, \ \ell\in{\mathbb N}$. 

\noindent
{\sc Case C.} \ The parameters are in the vertices of the boundary of
the domain of stability, i.e. \ $\alpha\beta\gamma =-1$ \ 
and \ $|\alpha|=|\beta|=|\gamma|=1$.

\begin{Thm}
   \label{main}
\ Let \ $\{\varepsilon_{k,\ell}:k,\ell\in{\mathbb N}\}$ \ be independent random
variables with \ ${\mathsf E}\,\varepsilon_{k,\ell}=0$ \ and \
$\Var\,\varepsilon_{k,\ell}=1$.
\ Assume that model \eqref{model} holds and for \ $n\in{\mathbb N}$ \
consider the piecewise constant random field
\begin{equation*}
Y^{(n)}(s,t):=X_{[ns],[nt]}, \qquad s,t\in{\mathbb R}, \quad s,t\geq 0.
\end{equation*} 

\noindent
If   \ $\ |\alpha |<1, \ |\beta |<1$ \ and \ $|\gamma|<1$, \
$|1+\alpha ^2-\beta ^2  -\gamma ^2|>2|\alpha +\beta\gamma|$ \ and \
$1-\beta ^2>|\alpha +\beta\gamma|$ \
then  
\begin{equation*}
\lim_{n\to\infty} \Var
\big(Y^{(n)}(s,t)\big)=\sigma^2_{\alpha,\beta,\gamma}:=
\big((1\!+\!\alpha\!+\!\beta\!-\!\gamma) 
(1\!+\!\alpha\!-\!\beta\!+\!\gamma)(1\!-\!\alpha\!+\!\beta\!+\!\gamma)
(1\!-\!\alpha\!-\!\beta\!-\!\gamma )\big)^{-1/2}.   
\end{equation*}

\noindent
If   \ $\ |\alpha |<1, \ |\beta |<1$ \ and \ $|\gamma|<1$ \ and in case
\ $\alpha\beta\gamma \geq 0$ \ equation  \ $|\alpha|+|\beta|+|\gamma|=1$,
\ while in case \ $\alpha\beta\gamma < 0$ \ equation  \
$|\alpha|+|\beta|-|\gamma|=1$ \ holds then  
\begin{equation*}
\lim_{n\to\infty} \frac 1{n^{1/2}}\Var \big(Y^{(n)}(s,t)\big)=\frac
{\big((1-|\alpha|)s\big)^{1/2}\land
  \big((1-|\beta|)t\big)^{1/2}}{\pi^{1/2}(|\alpha|+|\beta|)^{1/2}
  (1-|\alpha|)(1-|\beta|)}.   
\end{equation*}

\noindent
If \ $\alpha\beta\gamma \leq 0$ \ and \ $|\alpha|=1, \ |\beta|=|\gamma|<1$ \
or \ \ $|\beta|=1, \ |\alpha|=|\gamma|<1$ \ then 

\begin{equation*}
\lim_{n\to\infty} \frac 1n\Var \big(Y^{(n)}(s,t)\big)=\frac
s{1-\gamma^2} \qquad \text{or} \qquad \lim_{n\to\infty} \frac
1n\Var \big(Y^{(n)}(s,t)\big)=\frac t{1-\gamma^2},
\end{equation*}
respectively.

\noindent
If \ $\alpha\beta\gamma =-1$ \ and \ $|\alpha|=|\beta|=|\gamma|=1$ \
then  
\begin{equation*}
\lim_{n\to\infty} \frac 1{n^2}\Var \big(Y^{(n)}(s,t)\big)=st.
\end{equation*}
\end{Thm}

Observe that in the last case the limit of the variances is exactly
the variance of the standard Wiener sheet. This result is quite
natural, as e.g. for \ $\alpha=\beta=-\gamma=1$ \ model equation \eqref{model}
reduces to \
$\Delta_{1,1}\Delta_{2,1}X_{k,\ell}=\varepsilon_{k,\ell}$. 

We remark that results given Theorem \ref{main} do not cover all
possible locations of the parameters on the boundary of the domain of
stability. Some results on the missing cases 
can be found in Section \ref{sec:sec4}. 

The aim of the following discussion is to show that it suffices to
prove Theorem \ref{main} for \ $\alpha \geq 0$, \ $\beta \geq 0$ \ and
\ $\gamma \geq 0$ \ if \ $\alpha\beta\gamma\geq 0$ \ and for  \
$\alpha >0$, \ $\beta >0$ \ and \ $\gamma < 0$ \ if \
$\alpha\beta\gamma <0$. \ First we note that direct calculations imply
\begin{align}\label{marep}
  X_{k,\ell}&=\sum_{i=1}^k\sum_{j=1}^{\ell}G(k-i,\ell-j;\alpha,\beta,\gamma)
            \varepsilon_{i,j}\\ 
           &=\sum_{i=1}^k\sum_{j=1}^{\ell}\binom{k+\ell -i-j}{\ell
             -j}\alpha^{k-i}\beta^{\ell-j}F\bigg(i-k,j-\ell\,;i+j-k-\ell\,;
           -\frac{\gamma}{\alpha\beta}\bigg) \varepsilon_{i,j}, \label{marepf}
 \end{align}
$k,\ell \geq 1$, \ where \eqref{marepf} holds only for \
$\alpha\beta\ne 0$,
\begin{equation*}
G(m,n;\alpha,\beta,\gamma):=\sum_{r=0}^{m\land n} \frac{(m+n-r)!}{(m-r)!(n-r)!r!}
                     \alpha^{m-r}\beta^{n-r}\gamma ^r, \qquad m,n\in
                     {\mathbb N}\cup \{0\},
\end{equation*}
and \ $F(-n,b;c;z)$ \ is the Gauss hypergeometric
function defined by
\begin{equation*}
F(-n,b;c;z):=\sum_{r=0}^{n }
\frac {(-n)_r(b)_r}{(c)_rr!}z^r, \qquad n\in{\mathbb N},  \quad
b,c,z\in{\mathbb C},
\end{equation*}
and \ $(a)_r:=a(a+1)\dots (a+r-1)$ \ (for the definition in more
general cases see e.g. \citet{be}).

Observe that as for \ $m,n\in{\mathbb N}$ \ we have \
$F(-n,-m;-n-m;1)=\binom{m+n}n^{-1}$ \ and \ $F(-n,-m;-n-m;0)=1$, \
moving average representations of the doubly geometric model of
\citet{martin1} and of the spatial models studied by \citet{paul} and
\citet{bpz1,bpz2}, respectively, are special forms of \eqref{marepf}.

Now, put \
$\widetilde\varepsilon_{k,\ell}:=(-1)^{k+\ell}\varepsilon_{k,\ell}$ \
for \ $k,\ell \in{\mathbb N}$. \ Then \
$\{\widetilde\varepsilon_{k,\ell}:k,\ell\in{\mathbb N}\}$ \ are
independent random variables with \ ${\mathsf
  E}\,\widetilde\varepsilon_{k,\ell}=0$ \ and \
$\Var\,\widetilde\varepsilon_{k,\ell}=1$.
 \ Consider the process \ $\{\widetilde
X_{k,\ell}:k,\ell\in{\mathbb Z},\,k, \ell\geq0\}$ \  defined as 
 \begin{equation*}
  \begin{cases} \widetilde X_{k,\ell}=
                     -\alpha \widetilde X_{k-1,\ell}-\beta \widetilde
                     X_{k,\ell-1}+\gamma \widetilde
                     X_{k-1,\ell-1}+\widetilde\varepsilon_{k,\ell},  
                    & \text{for \ $k,\ell \geq 1$,}\\
                    \widetilde X_{k,0}=\widetilde X_{0,\ell}=0,  &
                    \text{for \ $k,\ell\geq0$.} 
                   \end{cases}
 \end{equation*}
Then by representation \eqref{marep} for \ $k,\ell\in{\mathbb N}$ \
we have
\begin{align*}
\widetilde X_{k,\ell}&=\sum_{i=1}^k\sum_{j=1}^{\ell}
  G(k-i,\ell-j;-\alpha,-\beta,\gamma)
  \widetilde\varepsilon_{i,j} \\
&=\sum_{i=1}^k\sum_{j=1}^{\ell}
  (-1)^{k+\ell-i-j}G(k-i,\ell-j;\alpha,\beta,\gamma)
  \widetilde\varepsilon_{i,j}=(-1)^{k+\ell}X_{k,\ell},
\end{align*}
hence \ $\Var \widetilde X_{k,\ell}=\Var X_{k,\ell}$. 

Next, put \
$\widehat\varepsilon_{k,\ell}:=(-1)^k\varepsilon_{k,\ell}$ \
for \ $k,\ell \in{\mathbb N}$. \ Then \
$\{\widehat\varepsilon_{k,\ell}:k,\ell\in{\mathbb N}\}$ \ are again
independent random variables with \ ${\mathsf
  E}\,\widehat\varepsilon_{k,\ell}=0$ \ and \ 
$\Var\,\widehat\varepsilon_{k,\ell}=1$.
 \ Consider the process \ $\{\widehat
X_{k,\ell}:k,\ell\in{\mathbb Z},\,k, \ell\geq0\}$ \  defined as 
 \begin{equation*}
  \begin{cases} \widehat X_{k,\ell}=
                     -\alpha \widehat X_{k-1,\ell}+\beta \widehat
                     X_{k,\ell-1}-\gamma \widehat
                     X_{k-1,\ell-1}+\widehat\varepsilon_{k,\ell},  
                    & \text{for \ $k,\ell \geq 1$,}\\
                    \widehat X_{k,0}=\widehat X_{0,\ell}=0,  &
                    \text{for \ $k,\ell\geq0$.} 
                   \end{cases}
 \end{equation*}
Then by representation \eqref{marep} for \ $k,\ell\in{\mathbb N}$ \
we have
\begin{align*}
\widehat X_{k,\ell}&=\sum_{i=1}^k\sum_{j=1}^{\ell}
 G(k-i,\ell-j;-\alpha,\beta,-\gamma)
\widehat\varepsilon_{i,j} \\
&=\sum_{i=1}^k\sum_{j=1}^{\ell} (-1)^{k-i}
 G(k-i,\ell-j;\alpha,\beta,\gamma)
\widehat\varepsilon_{i,j}=(-1)^kX_{k,\ell}, 
\end{align*}
hence \ $\Var \widehat X_{k,\ell}=\Var X_{k,\ell}$. 

In a similar way we have \ $\bar X_{k,\ell}=(-1)^{\ell}X_{k,\ell}$, \
so \ $\Var \bar X_{k,\ell}=\Var X_{k,\ell}$, \
where \ $\{\bar X_{k,\ell}:k,\ell\in{\mathbb Z},\,k, \ell\geq0\}$ \
is defined as 
 \begin{equation*}
  \begin{cases} \bar X_{k,\ell}=
                     \alpha \bar X_{k-1,\ell}-\beta \bar
                     X_{k,\ell-1}-\gamma \bar
                     X_{k-1,\ell-1}+\bar\varepsilon_{k,\ell},  
                    & \text{for \ $k,\ell \geq 1$,}\\
                    \bar X_{k,0}=\bar X_{0,\ell}=1.  &
                    \text{for \ $k,\ell\geq0$,} 
                   \end{cases}
 \end{equation*}
with \ $\bar\varepsilon_{k,\ell}:=(-1)^{\ell}\varepsilon_{k,\ell}$.

\section{Upper bounds for the covariances}
  \label{sec:sec2}
\markboth{\rm Baran}{\rm On the variances of a spatial
  unit root model} 

By representations \eqref{marep} and \eqref{marepf} we obtain that for 
\ $k_1,\ell_1,k_2,\ell_2\in {\mathbb N}$ \  and  \
$\alpha,\beta,\gamma \in{\mathbb R}$ \ we have
 \begin{align}\label{eq:eq2.1}
  \Cov \big(X_{k_1,\ell_1}&,X_{k_2,\ell_2})\!=\!\!\!\sum_{i=1}^{k_1\land k_2}
    \sum_{j=1}^{\ell_1\land
      \ell_2}\!G(k_1-i,\ell_1-j;\alpha,\beta,\gamma) 
     G(k_2-i,\ell_2-j;\alpha,\beta,\gamma)\\
           &=\!\!\!\sum_{i=1}^{k_1\land k_2}\sum_{j=1}^{\ell_1\land
             \ell_2} \binom{k_1+\ell_1 -i-j}{\ell_1 -j}
            \binom{k_2+\ell_2 -i-j}{\ell_2
              -j}\alpha^{k_1+k_2-2i}\beta^{\ell_1+\ell_2-2j} \label{eq:eq2.2}\\
          &\phantom{=\!\!\!\sum_{i=1}^{k_1\land
              k_2}}
           \times \!F\bigg(i\!-\!k_1,j\!-\!\ell_1;i\!+\!j\!-\!k_1\!-\!\ell_1;
           -\frac{\gamma}{\alpha\beta}\bigg)
           F\bigg(i\!-\!k_2,j\!-\!\ell_2;i\!+\!j\!-\!k_2\!-\!\ell_2; 
           -\frac{\gamma}{\alpha\beta}\bigg), \nonumber
 \end{align}  
where \ $x\land y:=\min\{x,y\}, \ x,y \in {\mathbb R}$, \ and
\eqref{eq:eq2.2} holds only for \ $\alpha\beta\ne 0$.

To obtain a more convenient form of the covariances we prove the
following Lemma.

\begin{Lem}
   \label{binoform}
\ Let \ $n,m$ \ be nonnegative integers and let \ $\alpha, \beta,
\gamma \in{\mathbb R}$ \ such that \ $0\leq \alpha, \beta <1 $ \
and \ $\alpha\beta+\gamma\geq 0$. \ Then
\begin{equation*}
G(m,n;\alpha,\beta,\gamma)
=\left (\frac{\alpha+\gamma}{1-\beta}\right )^m\!\!{\mathsf P}\left (
\xi_n^{(\beta)}+\eta_m^{\big(\frac{\alpha\beta+\gamma}{\alpha+\gamma}\big)}=n
\right) 
=\left (\frac{\beta+\gamma}{1-\alpha}\right )^n\!\!{\mathsf P}\left (
\xi_m^{(\alpha)}+\eta_n^{\big(\frac{\alpha\beta+\gamma}{\beta+\gamma}\big)}=m
\right), 
\end{equation*}
where \ $\xi_n^{(\nu)}$ \ and \ $\eta_m^{(\mu)},$ \  $0\leq \mu,\nu
\leq 1$, \ are independent binomial random variables with parameters \ 
$(n,\nu)$ \ and \ $(m,\mu)$, \ respectively, if \ $m,n\in{\mathbb N}$, \
and \ $\xi_0^{(\nu)}=\eta_0^{(\mu)}:=0$.
\end{Lem}

\noindent
{\bf Proof} \ As in cases \ $\alpha\beta\gamma=0$ \ or  \
$\alpha\beta+\gamma=0$ \ the statement of the 
Lemma holds trivially, we may assume \ $\alpha\beta\gamma\ne 0$ \ and
\ $\alpha\beta+\gamma>0$.  \ Let \
$0<n\leq m$.  
\begin{align*}
G(m,n;&\,\alpha,\beta,\gamma)=
\sum_{r=0}^{m\land n} \frac{(m+n-r)!}{(m-r)!(n-r)!r!}
                     \alpha^{m-r}\beta^{n-r}\gamma ^r=\frac
                     {\big(\beta+|\gamma|\big)^n\big(\Sign(\gamma)\big)^m}{n!}
                      \\ &\times\sum_{r=0}^n \frac 
                     {(n+m-r)!}{(m-r)!}{\alpha_{\gamma}}
                     ^{m-r} {\mathsf P}\Big 
                     (\xi_n^{(\nu)}=n-r \Big) 
                     =\frac{\big(\beta+|\gamma|\big)^n
                       \big(\Sign(\gamma)\big)^m}{n!}  
                   \frac {\mathrm d^n\big(\alpha_{\gamma}^m
                     {\mathcal G}_n^{(\nu)}(\alpha_{\gamma})\big )}{\mathrm
                     d\alpha_{\gamma}^n},  
\end{align*} 
where \ $\nu:=\frac{\beta}{\beta+|\gamma|}$, \
$\alpha_{\gamma}:=\Sign(\gamma)\,\alpha$ \ and \ 
${\mathcal G}_n^{(\nu)}(x):=\big (\nu x +(1-\nu)\big )^n$ \ is the
generating function of \ $\xi_n^{(\nu)}$. \ From the other hand
\begin{equation*}
 \frac {\mathrm d^n\big(\alpha_{\gamma} ^m
                     {\mathcal G}_n^{(\nu)}(\alpha_{\gamma})\big )}{\mathrm d
                     \alpha_{\gamma} 
                     ^n}=n! \sum _{r=0}^n\binom nr \binom m{n-r}\nu
                   ^r\big (\nu\alpha_{\gamma} +(1-\nu)\big
                   )^{n-r}\alpha_{\gamma}^{m-(n-r)}, 
\end{equation*}
and as \ $\alpha\beta+\gamma<\alpha+\gamma$ \ we obtain 
\begin{align}
 \label{eq:eq2.3}
G&(m,n;\alpha,\beta,\gamma)=\sum
                     _{r=0}^n\binom nr \binom m{n-r}\beta 
                   ^r(\alpha\beta +\gamma)^{n-r}\alpha ^{m-n+r}\\
                   &=\left
                     (\frac{\alpha+\gamma}{1-\beta}\right )^m \sum
                   _{r=0}^n\binom nr \beta ^r(1-\beta )^{n-r}\binom
                   m{n-r}\left(\frac{\alpha\beta+\gamma}{\alpha+\gamma}\right
                 )^{n-r}
                 \left(\frac{\alpha(1-\beta)}{\alpha+\gamma}\right
                 )^{m-n+r} \nonumber  \\
                  &=\left (\frac{\alpha+\gamma}{1-\beta}\right )^m \sum
                   _{r=0}^n{\mathsf P}\Big (\xi_n^{(\beta)}=r\Big)
                   {\mathsf P}\left
                 (\eta_m^{\big(\frac{\alpha\beta+\gamma}{\alpha+\gamma}\big)}
                 =n\!-\!r   
                   \right)\!=\!\left (\frac{\alpha+\gamma}{1-\beta}\right
                   )^m\!\!{\mathsf P}\left(\xi_n^{(\beta)}\!+
                  \eta_m^{\big(\frac{\alpha\beta+\gamma}{\alpha+\gamma}\big)}\!=n
                   \right)\!. \nonumber
\end{align}
Moreover, 
\begin{align*}
 \sum_{r=0}^n&\binom nr \binom m{n-r}\beta^r(\alpha\beta
 +\gamma)^{n-r}\alpha ^{m-n+r}= \sum_{r=m-n}^m\binom mr \binom n{m-r}\alpha 
                   ^r(\alpha\beta +\gamma)^{m-r}\beta ^{n-m+r}\\
  &=\left (\frac{\beta+\gamma}{1-\alpha}\right )^n \!\!\sum
                   _{r=m-n}^m\!\!{\mathsf P}\Big (\xi_m^{(\alpha)}=r\Big)
                   {\mathsf P}\left
                 (\eta_n^{\big(\frac{\alpha\beta+\gamma}{\beta+\gamma}\big)}\!\!
                 =m\!-\!r   
                   \right)\!=\!\left (\frac{\beta+\gamma}{1-\alpha}\right
                   )^n\!\!{\mathsf P}\left(\xi_m^{(\alpha)}\!+
                  \eta_n^{\big(\frac{\alpha\beta+\gamma}{\beta+\gamma}\big)}\!=m
                   \right)\! \nonumber 
\end{align*}
that together with \eqref{eq:eq2.3} implies the statement of the
lemma. Case \ $n>m$ \ can be handled in a similar way. \proofend

\begin{Cor}
   \label{cor:cor1} \
If \ $0\leq\alpha,\beta<1$ \ and \ $\alpha+\beta+\gamma=1$ \ then
\begin{equation*}
G(m,n;\alpha,\beta,\gamma)
={\mathsf P}\Big (\xi_m^{(\alpha)}+\eta_n^{(1-\beta)}=m \Big)=
{\mathsf P}\Big (\xi_n^{(\beta)}+\eta_m^{(1-\alpha )}=n \Big).
\end{equation*}
\end{Cor}

The following lemma is a natural generalization of Theorem 2.4 of
\citet{bpz2}.

\begin{Lem}
   \label{binoappr}
\ Let \ $k,\ell\in {\mathbb N}$, \ let \ $0<\mu, \ \nu<1$ \ be  real numbers
and let \ $\xi_k^{(\nu)}$ \ and  \
$\eta_{\ell }^{(\mu)}$ \ be independent binomial random variables
with parameters \ $(k,\nu)$ \ and \ $(\ell,\mu)$, \
respectively. Further, let \ $S_{k,\ell}:=\xi_k^{(\nu)}+\eta_{\ell
}^{(\mu)}$ \ and let
\begin{equation*}
m_{k,\ell}:={\mathsf E} S_{k,\ell} ,\qquad
b_{k,\ell}:=\Var S_{k,\ell},\qquad  
    x_{j,k,\ell}:=(j-m_{k,\ell})/\sqrt{b_{k,\ell}}.
\end{equation*}
Then for all \ $k,\ell \in{\mathbb N}$ \ and \
$j\in\{0,1,\ldots,k+\ell\}$, \ we have
\begin{equation*}
\left|{\mathsf P}(S_{k,\ell}=j)
 -\frac1{\sqrt {2\pi b_{k,\ell}}}\exp\left (-x_{j,k,\ell}^2/2\right )\right|
    \leq\frac{C(\mu,\nu)}{b_{k,\ell}},
\end{equation*}
where \ $C(\mu,\nu)>0$ \ is a constant depending only on \
$\mu$ \ and \ $\nu$ \ (and not 
depending on  \ $k,\ell,j$).
\end{Lem}

\begin{Thm} \label{covbound}
\ If \ $|\alpha|+|\beta|+|\gamma|<1$ \ then 
\begin{equation*}
\big |\Cov (X_{k_1,\ell_1}, X_{k_2,\ell_2})\big |\leq \frac {\big (
  |\alpha|+|\beta|+|\gamma|)^{
    (|k_1-k_2|+|\ell_1-\ell_2|)/2}}{\big(1-
  (|\alpha|+|\beta|+|\gamma|)\big)^2}.   
\end{equation*}

If   \ $\ |\alpha |<1, \ |\beta |<1$ \ and \ $|\gamma|<1$ \ and in case
\ $\alpha\beta\gamma \geq 0$ \ equation  \ $|\alpha|+|\beta|+|\gamma|=1$,
\ while in case \ $\alpha\beta\gamma < 0$ \ equation  \
$|\alpha|+|\beta|-|\gamma|=1$ \ holds then  
\begin{equation*}
\big |\Cov (X_{k_1,\ell_1}, X_{k_2,\ell_2})\big |\leq \
C\big(\alpha,\beta\big) \sqrt{k_1+\ell_1+k_2+\ell_2}
 \end{equation*}
with some constant \ $C\big(\alpha,\beta\big)>0$.

If \ $\alpha\beta\gamma \leq 0$ \ and \ $|\alpha|=1, \ |\beta|=|\gamma|<1$ \
or \ \ $|\beta|=1, \ |\alpha|=|\gamma|<1$ \ then 
\begin{equation*}
\big |\Cov (X_{k_1,\ell_1}, X_{k_2,\ell_2})\big |\leq \ (k_1\land
k_2) \frac {|\gamma|^{|\ell_1 -\ell_2|}}{1-\gamma^2} \quad \text{or}
\quad
\big |\Cov (X_{k_1,\ell_1}, X_{k_2,\ell_2})\big |\leq \ (\ell_1\land
\ell_2) \frac {|\gamma|^{|k_1 -k_2|}}{1-\gamma^2},
\end{equation*}
respectively.

If \ $\alpha\beta\gamma =-1$ \ and \ $|\alpha|=|\beta|=|\gamma|=1$ \
then  
\begin{equation*}
\Cov (X_{k_1,\ell_1}, X_{k_2,\ell_2})=(k_1\land k_2)
(\ell_1\land \ell_2)\,\alpha ^{|k_1-k_2|} \beta ^{|\ell_1-\ell_2|}.
\end{equation*}
\end{Thm}

\noindent
{\bf Proof.} \ Let  \ $|\alpha|+|\beta|+|\gamma|<1$. \ Lemma
\ref{binoform} and \eqref{eq:eq2.1} imply
\begin{align*}
\big |\Cov (X_{k_1,\ell_1}, X_{k_2,\ell_2})\big |&\leq
\sum_{i=1}^{k_1\land k_2}\sum_{j=1}^{\ell_1\land \ell _2} \left
  (\frac{|\alpha|+|\gamma|}{1-|\beta|}\right )^{(k_1+k_2)/2-i} \left
  (\frac{|\beta|+|\gamma|}{1-|\alpha|}\right )^{(\ell_1+\ell_2)/2-j}
\\ &\leq
\left(\frac{|\alpha|+|\gamma|}{1-|\beta|}\right )^{|k_1-k_2|/2} \left
  (\frac{|\beta|+|\gamma|}{1-|\alpha|}\right )^{|\ell_1-\ell_2|/2}
\frac {(1-|\beta|)(1-|\alpha|)}{\big(1- (|\alpha|+|\beta|+|\gamma|)\big)^2} . 
 \end{align*}
Hence, as 
\begin{equation*}
\frac{|\alpha|+|\gamma|}{1-|\beta|} \leq |\alpha|+|\beta|+|\gamma|
\quad \text{and} \quad \frac{|\beta|+|\gamma|}{1-|\alpha|} \leq
|\alpha|+|\beta|+|\gamma| 
\end{equation*} 
hold, we obtain the first statement of the theorem.

Now, let  \ $|\alpha|< 1, \ |\beta|< 1, \ |\gamma|< 1$ \ and assume
that 
in case \ $\alpha\beta\gamma \geq 0$ \ equation  
\ $|\alpha|+|\beta|+|\gamma|=1$,
\ while in case \ $\alpha\beta\gamma < 0$ \ equation  \
$|\alpha|+|\beta|-|\gamma|=1$ \ holds. From the arguments of the
Introduction follows that it suffices to consider the case \
$0\leq \alpha,\beta <1, \ |\gamma |<1$, \ and \ $\alpha+\beta+\gamma=1$.
 \ Corollary 
\ref{cor:cor1} and \eqref{eq:eq2.1} imply
\begin{equation*}
\big |\Cov (X_{k_1,\ell_1}, X_{k_2,\ell_2})\big |\leq
\sum_{i=1}^{k_1\land k_2}\sum_{j=1}^{\ell_1\land \ell _2} {\mathsf
  P}\Big (\xi_{k_1-i}^{(\alpha)}+\eta_{\ell_1-j}^{(1-\beta )}=k_1-i
\Big){\mathsf P}\Big (\xi_{k_2-i}^{(\alpha)}+\eta_{\ell_2-j}^{(1-\beta
    )}=k_2-i \Big).  
\end{equation*}
Assume first \
$k_1\leq k_2$ \ and \ $\ell_1\leq \ell_2$ \ or \ $k_1> k_2$ \ and \
$\ell_1>\ell_2$. \ In this case using the notations and results of
Lemma \ref{binoappr} we have 
\begin{align*} 
\big |\Cov (&X_{k_1,\ell_1}, X_{k_2,\ell_2})\big |\leq\!\!\!
\sum_{i=0}^{k_1\land k_2-1}\sum_{j=0}^{\ell_1\land \ell _2-1} \!\!\!{\mathsf
  P}\Big (\xi_{|k_1-k_2|+i}^{(\alpha)}\!+\!
\eta_{|\ell_1-\ell_2|+j}^{(1-\beta )} =|k_1\!-\!k_2|\!+\!i\Big){\mathsf P} 
\Big (\xi_i^{(\alpha)}\!+\!\eta_j^{(1-\beta )}\!=\!i \Big) \\
\leq& \sum_{i=2}^{k_1\land k_2-1}\sum_{j=2}^{\ell_1\land \ell _2-1}\frac
 1{2\pi\sqrt {b_{|k_1-k_2|+i,|\ell_1-\ell_2|+j}}\sqrt {b_{i,j}}} \exp
 \Big (- \frac
 {x_{|k_1-k_2|+i, |k_1-k_2|+i, |\ell_1-\ell_2|+j}^2}2-\frac
 {x_{i,i,j}^2}2\Big), \\ 
 &+C\big(\alpha,\beta\big) \Big (\sum_{i=2}^{k_1\land
  k_2-1} \frac 1i +\sum_{j=2}^{\ell _1\land  \ell _2-1} \frac
1j+\sum_{i=2}^{k_1\land k_2-1}\sum_{j=2}^{\ell_1\land \ell _2-1}
 b_{i,j}^{-3/2}\Big)+4, 
\end{align*}
where \ $C\big(\alpha,\beta\big)$ \ is a positive
constant and 
\begin{equation*}
b_{k,\ell}:=\alpha(1-\alpha)k+ \beta(1-\beta)\ell  \qquad \text{and} \qquad
x_{k,k,\ell}:= a_{k,\ell}/\sqrt{b_{k,\ell}} \quad \text{with} \quad
a_{k,\ell}:=(1-\alpha)k - (1-\beta)\ell. 
\end{equation*}
 Obviously,
\begin{equation*}
\sum_{i=2}^{k_1\land k_2-1} \frac 1i \leq 2\sqrt{k_1\land k_2}\leq
2\sqrt{k_1+ k_2}  \qquad 
\text{and} \qquad
\sum_{j=2}^{\ell _1\land  \ell _2-1} \frac 1j \leq 2\sqrt{\ell _1\land
  \ell _2}\leq 2\sqrt{\ell _1+ \ell _2}.
\end{equation*}
Further, we have
\begin{align*}
\sum_{i=2}^{k_1\land k_2-1}&\sum_{j=2}^{\ell_1\land \ell _2-1}
 b_{i,j}^{-3/2}=\sum_{i=2}^{k_1\land k_2-1}\sum_{j=2}^{\ell_1\land
   \ell _2-1} \big(\alpha(1-\alpha)i +\beta(1-\beta)j
 \big)^{-3/2} \\
 &\leq\!\!\! \sum_{i=2}^{k_1\land k_2-1}\!\! \int\limits_1^{\ell_1\land
   \ell _2-1}\!\!\!\!\!\! \big(\alpha(1-\alpha)i +\beta(1-\beta)t
 \big)^{-3/2}{\mathrm d}t \leq \frac
 2{\beta(1-\beta)}\!\!\!\!\sum_{i=2}^{k_1\land k_2-1}\!\!\!\!
 \big(\alpha(1-\alpha)i\big)^{-1/2} \\
 &\leq\frac 2{\beta(1-\beta)} \int\limits_1^{k_1\land k_2-1}
 \big(\alpha(1-\alpha)s\big)^{-1/2}{\mathrm d}s\leq \frac 
 {4\big(\alpha(1-\alpha)\big)^{1/2}}{\alpha\beta(1-\alpha)(1-\beta)}
 \sqrt{k_1\land k_2 }.
\end{align*}
Hence,
\begin{equation}
   \label{eq:eq2.4} 
\big |\Cov (X_{k_1,\ell_1}, X_{k_2,\ell_2})\big |\leq
C\big(\alpha,\beta\big) \sqrt{k_1+\ell _1+k_2+ \ell
  _2} +4H_{\alpha,\beta}(k_1,\ell_1,k_2,\ell_2),
\end{equation}
with
\begin{align*}
H_{\alpha,\beta}&(k_1,\ell_1,k_2,\ell_2) \\ :=&\int\limits_1^{k_1\land
  k_2}\int\limits_1^{\ell_1\land \ell_2} 
\frac 1{2\pi\sqrt {b_{|k_1-k_2|+s,|\ell_1-\ell_2|+t}}\sqrt {b_{s,t}}} \exp
 \Big (- \frac
 {x_{|k_1-k_2|+s, |k_1-k_2|+s,|\ell_1-\ell_2|+t}^2}2-\frac
 {x_{s,s,t}^2}2\Big){\mathrm dt}{\mathrm d}s. 
\end{align*}
It is easy to see that 
\begin{align*}
&H_{\alpha,\beta}(k_1,\ell_1,k_2,\ell_2) \leq \frac
1{(\alpha+\beta)(1-\alpha)(1-\beta)} \\ &\phantom{=}\times
\!\!\!\!\!\!\!\! 
\!\int\limits_{b_{1,1} }^{b_{k_1\land k_2,\ell_1\land
  \ell_2}}\int\limits_{a_{1,\ell_1\land
  \ell_2}}^{a_{k_1\land k_2,1}} \!\!\!\!\!\frac
1{2\pi\sqrt{(b_{|k_1-k_2|,|\ell_1-\ell_2|}+u)u}}  \exp\!\bigg(\!\!
    -\frac{(a_{|k_1-k_2|,|\ell_1-\ell_2|}+v)^2}{2(b_{|k_1-k_2|,|\ell_1-\ell_2|}+u)}
    -\frac{v^2}{2u}  \bigg) {\mathrm d}v {\mathrm d}u. 
\end{align*}
Now, for some real constants \ $a<b$ \ and \ $q, \varrho$ \ we have
\begin{align*}
\int\limits_a^b \!\!
  \exp\bigg(\!\!
    -\frac{(\varrho+v)^2}{2(q+u)}
    -\frac{v^2}{2u}  \bigg) {\mathrm d}v =&
   \frac {\sqrt{\pi(q+u)u}}{\sqrt{2(q+2u)}}  \exp\bigg(
    -\frac{\varrho^2}{2(q+2u)}\bigg)
    \\ &\times \left (\widetilde\Phi \bigg (\frac {(q+2u)b+
        \varrho u}{\sqrt{2(q+2u)
          (q+u)u}}\bigg)\!-\!
      \widetilde\Phi \bigg (\frac {(q+2u)a+
        \varrho u}{\sqrt{2(q+2u)
          (q+u)u}}\bigg) \right ),  
\end{align*}
where  \ $\widetilde\Phi(x):=2\Phi (\sqrt {2}x)-1, \ x\in{\mathbb
  R},$ \ is the Gauss error function defined with the help of the cdf \
$\Phi (x)$ \ of the standard normal distribution. Hence 
\begin{align*}
H_{\alpha,\beta}(k_1,\ell_1,k_2,\ell_2) \leq&  \frac
2{\sqrt{2\pi}(\alpha+\beta)(1-\alpha)(1-\beta)} \\
&\times\int\limits_{b_{1,1} }^{b_{k_1\land k_2,\ell_1\land
  \ell_2}}  \frac 1{\sqrt{b_{|k_1-k_2|,|\ell_1-\ell_2|}\!+\!2u}}  \exp\bigg(
    -\frac{a_{|k_1-k_2|,|\ell_1-\ell_2|}^2}{2(b_{|k_1-k_2|,|\ell_1-\ell_2|}\!+\!2u)}\bigg)
    {\mathrm d}u \\
\leq & \frac
2{\sqrt{2\pi}(\alpha+\beta)(1-\alpha)(1-\beta)}
 \int\limits_{b_{1,1} }^{b_{k_1\land k_2,\ell_1\land
  \ell_2}}  \frac 1{\sqrt{b_{|k_1-k_2|,|\ell_1-\ell_2|}+2u}}{\mathrm d}u \\
\leq & \frac{\sqrt{
  2\alpha(1-\alpha)(k_1+k_2)+2\beta(1-\beta)(\ell_1+\ell_2)}}
{\sqrt{\pi}(\alpha+\beta)(1-\alpha)(1-\beta)} \!
\leq \!\! \frac{\sqrt{k_1+\ell_1+k_2+\ell_2}}
{\sqrt{2\pi}(\alpha+\beta)(1-\alpha)(1-\beta)}
\end{align*}
that together with \eqref{eq:eq2.4} implies the second statement of
the theorem.  Cases $k_1\leq k_2, \ \ell_1> \ell_2$ \ and \ $k_1> k_2,
\ \ell_1\leq \ell_2$ \ can be handled in a similar way.

Further, let \ $\alpha\beta\gamma <0$ \ and \ $|\alpha|=1, \
|\beta|=|\gamma|<1$ \ or \ $|\beta|=1, \
|\alpha|=|\gamma|<1$. \ In this case \ $-\gamma/(\alpha\beta)=1$. \ As
for \ $n,m\in{\mathbb N}$ \ one has \
$F(-n,-m;-n-m;1)=\binom{m+n}n^{-1}$, \ representation \eqref{eq:eq2.2} implies
\begin{equation}
   \label{eq:eq2.5}
\Cov (X_{k_1,\ell_1}, X_{k_2,\ell_2})=\sum_{i=1}^{k_1\land
  k_2}\sum_{j=1}^{\ell_1\land 
             \ell_2} \alpha^{k_1+k_2-2i}\beta^{\ell_1+\ell_2-2j}.
\end{equation}
Obviously, \eqref{eq:eq2.5} also holds if \ $|\alpha|=1, \
\beta=\gamma=0$ \ or \ $|\beta|=1, \
\alpha=\gamma=0$. \
Hence, e.g. if \ $|\alpha|=1, \ |\beta|=|\gamma|<1$
\begin{equation*}
\big | \Cov (X_{k_1,\ell_1}, X_{k_2,\ell_2}) \big | \leq (k_1\land
k_2)\!\!\sum_{j=1}^{\ell_1\land 
             \ell_2} |\gamma |^{\ell_1+\ell_2-2j} \leq (k_1\land
k_2)|\gamma |^{|\ell_1-\ell_2|}\!\!\sum_{j=0}^{\ell_1\land 
             \ell_2-1}\!\!\!\! \gamma ^{2j}\leq (k_1\land k_2)
           \frac{|\gamma|^{|\ell_1-\ell_2|}}{1-\gamma ^2}. 
\end{equation*}

Finally, if  \ $\alpha\beta\gamma =-1$ \ and \ $|\alpha|=|\beta|=|\gamma|=1$ \
then \  $-\gamma/(\alpha\beta)=1$, \ so
\begin{equation}
   \label{eq:eq2.6}
\Cov (X_{k_1,\ell_1}, X_{k_2,\ell_2})=\sum_{i=1}^{k_1\land k_2}\sum_{j=1}^{\ell_1\land
             \ell_2} \alpha^{k_1+k_2}\beta^{\ell_1+\ell_2}=(k_1\land k_2)
(\ell_1\land \ell_2)\,\alpha ^{|k_1-k_2|} \beta ^{|\ell_1-\ell_2|}
\end{equation}
that completes the proof. \proofend

\markboth{\rm Baran}{\rm On the variances of a spatial
  unit root model}

\section{Proof of Theorem \ref{main}}
  \label{sec:sec3}
\markboth{\rm Baran}{\rm On the variances of a spatial
  unit root model} 
 
According to the results of the Introduction we may assume \
$\alpha\geq 0, \ \beta\geq0$ \ and \ $\gamma \geq 0$ \ if \
$\alpha\beta\gamma \geq 0$ \ and  \ $\alpha> 0, \ \beta>0$ \ and \
$\gamma <0$, \ otherwise. 

Let  \ $0\leq\alpha ,\beta <1$ \ and \ $|\gamma|<1$, \
$|1+\alpha ^2-\beta ^2  -\gamma ^2|>2|\alpha +\beta\gamma|$ \ and \
$1-\beta ^2>|\alpha +\beta\gamma|$. \ Representation \eqref{eq:eq2.1}
directly implies 
\begin{equation}
   \label{eq:eq3.1}
\lim_{n\to\infty}\Var \big(Y^{(n)}(s,t)\big)=\sum_{i=0}^{\infty}
    \sum_{j=0}^{\infty}\Big (G(i,j;\alpha,\beta,\gamma)\Big)^2.
\end{equation} 
To show that the right hand side of \eqref{eq:eq3.1} equals \ $\sigma
^2_{\alpha,\beta,\gamma}$ \ consider the stationary solution \
$X_{k,\ell}^*$ \ of the equation
\begin{equation*}
X_{k,\ell}^*=\alpha X_{k-1,\ell}^*+\beta X_{k,\ell -1}^*+\gamma
X_{k-1,\ell-1}^*+\varepsilon_{k,\ell}^*, \quad k,\ell \in {\mathbb Z},
\end{equation*} 
where \ $\{\varepsilon^*_{k,\ell}:k,\ell\in{\mathbb Z}\}$ \ are
independent random 
variables with \ ${\mathsf E}\,\varepsilon^*_{k,\ell}=0$ \ and \
$\Var\,\varepsilon^*_{k,\ell}=1$. \ 
As the model is stable, \ $X_{k,\ell}^*$ \ has the following ${\mathsf
  L}^2$-convergent infinite moving average representation (see
\citet[Lemma 5.1]{tj2}) 
\begin{equation*}
X_{k,\ell}^*=\sum_{i=0}^{\infty}
    \sum_{j=0}^{\infty}G(i,j ;\alpha,\beta,\gamma)
                     \varepsilon _{k-i,\ell-j}.
\end{equation*}
Hence, 
\begin{equation*}
\Var \big(X_{k,\ell}^*\big)=\sum_{i=0}^{\infty}
    \sum_{j=0}^{\infty}\Big (G(i,j;\alpha,\beta,\gamma)\Big)^2.
\end{equation*} 
On the other hand, using the results of \citet{br2} one can easily show that
\ $\Var \big(X_{k,\ell}\big)=\sigma^2_{\alpha,\beta,\gamma}$.

Further, let  \ $0\leq\alpha ,\beta <1$, \ $|\gamma|<1$, \ and \
$\alpha+\beta+\gamma=1$. \ 
Corollary \ref{cor:cor1} and \eqref{eq:eq2.1} imply 
\begin{equation*}
\Var \big(Y^{(n)}(s,t)\big)=\sum _{k=0}^{[ns]-1}\sum
_{\ell=0}^{[nt]-1} {\mathsf P}^2\Big (\xi_k^{(\alpha)}+\eta_{\ell
}^{(1-\beta)}=k \Big).
\end{equation*} 
Hence, to find the limit on \ $n^{-1/2}\Var \big(Y^{(n)}(s,t)\big)$ \
as \ $n\to\infty$, \ one can use the local version of the central
limit theorem given in
Lemma \ref{binoappr} that yields approximation 
\begin{align*}
\Var \big(Y^{(n)}(s,t)\big)\approx \widetilde E^{(n)}_{\alpha,\beta}(s,t):=&\sum
_{k=1}^{[ns]-1}\sum_{\ell=1}^{[nt]-1} \frac 1{2\pi b_{k,\ell}}
\exp \big (- x_{k,k,\ell}^2\big) \\
=&\int\limits_1^{[ns]}\int\limits_1^{[nt]}
\frac 1{2\pi b_{[y],[z]}} \exp \big (- x_{[y],[y],[z]}^2\big){\mathrm
  d}z{\mathrm d}y.
\end{align*} 
Direct calculations show that for the error \
\begin{equation*}
\widetilde\Delta^{(n)}_{\alpha,\beta}(s,t):=\Var
\big(Y^{(n)}(s,t)\big)-\widetilde E^{(n)}_{\alpha,\beta}(s,t)
\end{equation*}
of the approximation we have
\begin{equation}
     \label{eq:eq3.2}
\big|\widetilde\Delta^{(n)}_{\alpha,\beta}(s,t)\big|\leq  C(\alpha,\beta)
\left (1+\sum _{k=2}^{[ns]-1}\sum_{\ell=2}^{[nt]-1} \frac 1{b_{k,\ell}^2} + \sum
_{k=2}^{[ns]-1}\sum_{\ell=2}^{[nt]-1} \frac 1{\sqrt{2\pi}
b_{k,\ell}^{3/2}} \exp \bigg (-\frac{x_{k,k,\ell}^2}2\bigg)\right), 
\end{equation}
where \ $C(\alpha,\beta)$ \ is a positive constant. 
Now, similarly to the proof of Theorem \ref{covbound} one can verify that
\begin{equation}
     \label{eq:eq3.3}
\sum_{k=2}^{[ns]-1}\sum_{\ell=2}^{[nt]-1} \frac 1{b_{k,\ell}^2} \leq \frac {\ln
  \big(\alpha(1-\alpha)([ns]-1) 
  +\beta(1-\beta)\big)}{\alpha\beta(1-\alpha)(1-\beta)}
 \leq \frac {\ln ([ns])}{\alpha\beta(1-\alpha)(1-\beta)} .  
\end{equation}
Further, 
\begin{align}
   \label{eq:eq3.4}
\sum_{k=2}^{[ns]-1}\sum_{\ell=2}^{[nt]-1} &\frac 1{b_{k,\ell}^{3/2}}
\exp \bigg (- \frac{x_{k,k,\ell}^2}2\bigg) \leq 4\int\limits_1^{[ns]}
\int\limits_1^{[nt]} \frac 1{b_{y,z}^{3/2}}
\exp \bigg (- \frac{x_{y,y,z}^2}2\bigg){\mathrm d}z{\mathrm d}y \\
&\leq \frac 4{(\alpha +\beta)(1-\alpha )(1-\beta )}
\int\limits_{b_{1,1}}^{b_{[ns],[nt]}}  
\int\limits_{a_{1,[nt]}}^{a_{[ns],1}} \frac 1{u^{3/2}}
\exp \bigg (- \frac{v^2}{2u}\bigg){\mathrm d}v{\mathrm d}u. \nonumber
\end{align}
Again, for some real constants \ $a<b$ \ and \ $m>0$
\begin{equation}
   \label{eq:eq3.5}
\int\limits_a^b \exp \bigg (-
\frac{v^2}{mu}\bigg){\mathrm d}v=\frac {\sqrt{\pi mu}}2 
\left (\widetilde\Phi \bigg (\frac b{\sqrt{mu}} \bigg) -
\widetilde\Phi \bigg (\frac a{\sqrt{mu}} \bigg)\right)
\end{equation}
holds, so using \eqref{eq:eq3.4} and \eqref{eq:eq3.5} with \ $m=2$ \ we have
\begin{align*}
\sum_{k=2}^{[ns]-1}\sum_{\ell=2}^{[nt]-1} &\frac 1{b_{k,\ell}^{3/2}}
\exp \bigg (\!\! - \frac{x_{k,k,\ell}^2}2\bigg) \leq \frac
{4\sqrt{2\pi}}{(\alpha \!
  +\!\beta)(1\!-\!\alpha )(1\!-\!\beta )} \ln \bigg(\frac
{[ns]+[nt]}{b_{1,1}}\bigg)
\end{align*}
that together with \eqref{eq:eq3.2} and \eqref{eq:eq3.3} implies
\begin{equation*}
  \lim_{n\to\infty}\frac 1{n^{1/2}}\widetilde\Delta^{(n)}_{\alpha,\beta}(s,t)=0.
\end{equation*}
Hence, \ $n^{-1/2}\Var \big(Y^{(n)}(s,t)\big)$ \ and \
$n^{-1/2}\widetilde E^{(n)}_{\alpha,\beta}(s,t)$ \ have the same
limit as \ $n\to\infty$.

Now, let 
\begin{equation*}
\Delta^{(n)}_{\alpha,\beta}(s,t):=\widetilde
E^{(n)}_{\alpha,\beta}(s,t)- E^{(n)}_{\alpha,\beta}(s,t), 
\end{equation*}
where
\begin{equation*}
E^{(n)}_{\alpha,\beta}(s,t):=\int\limits_1^{[ns]}\int\limits_1^{[nt]}
\frac 1{2\pi b_{y,z}} \exp \big (- x_{y,y,z}^2\big){\mathrm
  d}z{\mathrm d}y. 
\end{equation*}
Obviously,
\begin{equation}
   \label{eq:eq3.6}
\Delta^{(n)}_{\alpha,\beta}(s,t)=\Delta^{(n,1)}_{\alpha,\beta}(s,t)
+\Delta^{(n,2)}_{\alpha,\beta}(s,t),  
\end{equation}
where
\begin{align*}
\Delta^{(n,1)}_{\alpha,\beta}(s,t)&:=\frac
1{2\pi}\int\limits_1^{[ns]}\int\limits_1^{[nt]} 
\bigg(\frac 1{b_{[y],[z]}} \exp \Big (-
\frac{a_{[y],[z]}^2}{b_{[y],[z]}}\Big)-\frac 1{b_{y,z}} \exp \Big (-
\frac{a_{[y],[z]}^2}{b_{y,z}}\Big)\bigg){\mathrm d}z{\mathrm d}y,\\
\Delta^{(n,2)}_{\alpha,\beta}(s,t)&:=\frac
1{2\pi}\int\limits_1^{[ns]}\int\limits_1^{[nt]} 
\bigg(\frac 1{b_{y,z}} \exp \Big (-
\frac{a_{[y],[z]}^2}{b_{y,z}}\Big)-\frac 1{b_{y,z}} \exp \Big (-
\frac{a_{y,z}^2}{b_{y,z}}\Big)\bigg){\mathrm d}z{\mathrm d}y.
\end{align*}
As \ $\big|z-[z]\big|<1, \ z\in{\mathbb R}$, \ and for \ $z\geq0$ \ we
have \ $z\exp(-z)\leq 1$, \ and \ $|1-\exp(-z)|\leq |z|$, \ while  for
\ $z\geq 1$, \ $[z]>z/2$ \ holds, after short straightforward
calculations (see also \eqref{eq:eq3.3}) we obtain 
\begin{equation}
  \label{eq:eq3.7}
\big|\Delta^{(n,1)}_{\alpha,\beta}(s,t)\big|\leq
\int\limits_1^{[ns]}\int\limits_1^{[nt]} \frac
{b_{1,1}}{\pi b_{y,z}^2}{\mathrm d}z{\mathrm d}y \leq 
\frac{\ln([ns]+1)}{2\pi\alpha\beta(1-\alpha)(1-\beta)}.
\end{equation}
Further, using similar ideas as in the proof of \eqref{eq:eq3.7} we have
\begin{align*}
\big|\Delta^{(n,2)}_{\alpha,\beta}(s,t)\big|\leq& \frac 1{2\pi}\!\!
\int\limits_1^{[ns]}\int\limits_1^{[nt]}
\frac {\big|a_{y,z}^2\!-\!a_{[y],[z]}^2\big|}{b_{y,z}^2}
\exp \Big (\!\!-\frac{a_{y,z}^2\land a_{[y],[z]}^2}{b_{y,z}}\Big) {\mathrm
  d}z{\mathrm d}y \leq \! \frac
1{2\pi}\!\!\int\limits_1^{[ns]}\int\limits_1^{[nt]} 
\!\frac {(2-\alpha-\beta)^2}{b_{y,z}^2}{\mathrm d}z{\mathrm d}y
\nonumber \\ 
&+ \frac {2-\alpha-\beta}{\pi}
\int\limits_1^{[ns]}\int\limits_1^{[nt]}
\frac {|a_{y,z}|\land |a_{[y],[z]}|}{b_{y,z}^2}
\exp \Big (-\frac{a_{y,z}^2\land a_{[y],[z]}^2}{b_{y,z}}\Big){\mathrm
  d}z{\mathrm d}y \\ 
\leq&  \frac
4{\pi}\int\limits_1^{[ns]}\int\limits_1^{[nt]} 
\!\frac 1{b_{y,z}^2}{\mathrm d}z{\mathrm d}y
+\frac 2{\pi} \int\limits_1^{[ns]}\int\limits_1^{[nt]} {\mathcal
  X}_{\{|a_{y,z}|\land |a_{[y],[z]}|\geq 1\} }
\frac 1{b_{y,z} \big(|a_{y,z}|\land |a_{[y],[z]}|\big)}{\mathrm
  d}z{\mathrm d}y \nonumber\\
\leq&  \frac{8\ln([ns]+1)}{\pi\alpha\beta(1-\alpha)(1-\beta)}
+\frac 2{\pi} \int\limits_1^{[ns]}\int\limits_1^{[nt]} {\mathcal
  X}_{\{|a_{y,z}|\geq 1\} }
\frac 1{b_{y,z} |a_{y,z}|}{\mathrm d}z{\mathrm d}y \nonumber\\
\leq&  \frac{8\ln([ns]+1)}{\pi\alpha\beta(1-\alpha)(1-\beta)}
+\frac 2{\pi(\alpha+\beta)(1-\alpha)(1-\beta)}\int\limits
_{b_{1,1}}^{b_{[ns],[nt]}} \frac 1u{\mathrm d}u
\int\limits_{a_{1,[nt]}}^{a_{[ns],1}} {\mathcal X}_{\{|v|\geq 1\}
}\frac 1{|v|}{\mathrm d}v \nonumber \\
\leq&  \frac{8\ln([ns]+1)}{\pi\alpha\beta(1-\alpha)(1-\beta)}
+\frac {4\ln ([ns]+[nt])}{\pi(\alpha+\beta)(1-\alpha)(1-\beta)}\ln
\bigg(\frac {[ns]+[nt]}{b_{1,1}}\bigg), \nonumber 
\end{align*}
where \ ${\mathcal X}_H$ \ denotes the indicator function of a set \
$H$, that together with \eqref{eq:eq3.6} and \eqref{eq:eq3.7} implies
\begin{equation*}
  \lim_{n\to\infty}\frac 1{n^{1/2}}\Delta^{(n)}_{\alpha,\beta}(s,t)=0.
\end{equation*}
Hence, \ $n^{-1/2}\Var \big(Y^{(n)}(s,t)\big)$ \ and \
$n^{-1/2} E^{(n)}_{\alpha,\beta}(s,t)$ \ have the same
limit as \ $n\to\infty$.

Now, consider first the case \ $\alpha(1-\alpha)s\leq \beta(1-\beta)t$
\ implying  \ $\alpha(1-\alpha)[ns]+\beta(1-\beta)\leq
\alpha(1-\alpha)+\beta(1-\beta)[nt]$, \ if \ $n$ \ is large enough. In
this case 
\begin{equation*}
 E^{(n)}_{\alpha,\beta}(s,t)= \frac
 1{(\alpha+\beta)(1-\alpha)(1-\beta)}
\big(E^{(n,1)}_{\alpha,\beta}(s,t)+ E^{(n,2)}_{\alpha,\beta}(s,t)+
E^{(n,3)}_{\alpha,\beta}(s,t)\big), 
\end{equation*}
where
\begin{align*}
E^{(n,1)}_{\alpha,\beta}(s,t)&:=\frac
1{2\pi}\int\limits_{b_{1,1}}^{b_{[ns],1}}
\int\limits_{-u/\beta+(\alpha+\beta)(1-\alpha)/\beta}^{u/\alpha
  -(\alpha+\beta)(1-\beta)/\alpha} 
\frac 1u \exp\bigg(-\frac {v^2}u\bigg){\mathrm d}v{\mathrm d}u, \\
E^{(n,2)}_{\alpha,\beta}(s,t)&:=\frac
1{2\pi}\int\limits_{b_{[ns],1}}^{b_{1,[nt]}}
\int\limits_{-u/\beta+(\alpha+\beta)(1-\alpha)/\beta}^{-u/\beta+(\alpha+\beta)
  (1-\alpha)[ns]/\beta}
\frac 1u \exp\bigg(-\frac {v^2}u\bigg){\mathrm d}v{\mathrm d}u, \\
E^{(n,3)}_{\alpha,\beta}(s,t)&:=\frac
1{2\pi}\int\limits_{b_{1,[nt]}}^{b_{[ns],[nt]}}
\int\limits_{u/\alpha-(\alpha+\beta)(1-\beta)[nt]/\alpha}^{-u/\beta+(\alpha+\beta)
  (1-\alpha)[ns]/\beta}
\frac 1u \exp\bigg(-\frac {v^2}u\bigg){\mathrm d}v{\mathrm d}u. 
\end{align*}
Using \eqref{eq:eq3.5} with \ $m=1$, \ as  \
$\widetilde\Phi(-x)=-\widetilde\Phi(x)$, \ we have
\begin{align*}
E^{(n,1)}_{\alpha,\beta}(s,t)=&\frac
1{2\sqrt{\pi}}\int\limits_{b_{1,1}}^{b_{[ns],1}} \frac
1{2\sqrt{u}} \!\left(\widetilde\Phi 
  \bigg( \frac{\sqrt{u}}{\alpha}\!-\!\frac
  {(\alpha\!+\!\beta)(1\!-\!\beta)}{\sqrt{u}\alpha }\bigg )+\widetilde\Phi
  \bigg( \frac{\sqrt{u}}{\beta}\!-\!\frac
  {(\alpha\!+\!\beta)(1\!-\!\alpha)}{\sqrt{u}\beta }\bigg
  )\right){\mathrm d}u \\ 
=&\frac 1{2\sqrt{\pi}}\int\limits_{\sqrt{b_{1,1}}}^{\sqrt{b_{[ns],1}}} 
\left(\widetilde\Phi \bigg( \frac{w}{\alpha}\!-\!\frac
  {(\alpha\!+\!\beta)(1\!-\!\beta)}{w\alpha }\bigg )+\widetilde\Phi
  \bigg( \frac{w}{\beta}\!-\!\frac
  {(\alpha\!+\!\beta)(1\!-\!\alpha)}{w\beta }\bigg
  )\right){\mathrm d}w, \\
E^{(n,2)}_{\alpha,\beta}(s,t)=&\frac
1{2\sqrt{\pi}}\int\limits_{\sqrt{b_{[ns],1}}}^{\sqrt{b_{1,[nt]}}}  
\left(\widetilde\Phi \bigg( -\frac{w}{\beta}\!+\!\frac
  {(\alpha\!+\!\beta)(1\!-\!\alpha)[ns]}{w\beta }\bigg )+\widetilde\Phi
  \bigg( \frac{w}{\beta}\!-\!\frac
  {(\alpha\!+\!\beta)(1\!-\!\alpha)}{w\beta }\bigg
  )\right){\mathrm d}w, \\
E^{(n,3)}_{\alpha,\beta}(s,t)=&\frac
1{2\sqrt{\pi}}\!\!\int\limits_{\sqrt{b_{1,[nt]}}}^{\sqrt{b_{[ns],[nt]}}} \!\!  
\left(\widetilde\Phi \bigg( -\frac{w}{\beta}\!+\!\frac
  {(\alpha\!+\!\beta)(1\!-\!\alpha)[ns]}{w\beta }\bigg )-\widetilde\Phi
  \bigg( \frac{w}{\alpha}\!-\!\frac
  {(\alpha\!+\!\beta)(1\!-\!\beta)[nt]}{w\alpha }\bigg
  )\right){\mathrm d}w. 
\end{align*}
Combining similar terms we obtain
\begin{equation}
   \label{eq:eq3.8}
 E^{(n)}_{\alpha,\beta}(s,t)= \frac
 1{2\sqrt{\pi}(\alpha\!+\!\beta)(1\!-\!\alpha)(1\!-\!\beta)}
\big(F^{(n,1)}_{\alpha,\beta}(s,t)\!+\! F^{(n,2)}_{\alpha,\beta}(s,t)\!+\!
F^{(n,3)}_{\alpha,\beta}(s,t)\!+\!F^{(n,4)}_{\alpha,\beta}(s,t)\big), 
\end{equation}
where
\begin{alignat*}{2}
F^{(n,1)}_{\alpha,\beta}(s,t)&:=\!\!\!\!\!\!
\int\limits_{\sqrt{b_{1,1}}}^{\sqrt{b_{[ns],1}}}\!\!\!\!\!  
\widetilde\Phi \bigg( \frac{w}{\alpha}\!-\!\frac
  {(\alpha\!+\!\beta)(1\!-\!\beta)}{w\alpha }\bigg ){\mathrm d}w,\ \
F^{(n,3)}_{\alpha,\beta}(s,t)&:=\!\!\!\!\!\!\!\!
\int\limits_{\sqrt{b_{[ns],1}}}^{\sqrt{b_{[ns],[nt]}}}\!\!\!\!\!\!   
\widetilde\Phi \bigg(
\frac{(\alpha\!+\!\beta)(1\!-\!\alpha)[ns]}{w\beta }
\!-\!\frac{w}{\beta}\bigg ){\mathrm d}w , \\ 
F^{(n,2)}_{\alpha,\beta}(s,t)&:=\!\!\!\!\!\!
\int\limits_{\sqrt{b_{1,1}}}^{\sqrt{b_{1,[nt]}}}\!\!\!\!\!
\widetilde\Phi \bigg( \frac{w}{\beta}\!-\!\frac
  {(\alpha\!+\!\beta)(1\!-\!\alpha)}{w\beta }\bigg ){\mathrm d}w, \ \
F^{(n,4)}_{\alpha,\beta}(s,t)&:=\!\!\!\!\!\!\!\!
\int\limits_{\sqrt{b_{1,[nt]}}}^{\sqrt{b_{[ns],[nt]}}}\!\!\!\!\!\!   
\widetilde\Phi \bigg(\frac{(\alpha\!+\!\beta)(1\!-\!\beta)[nt]}{w\alpha }
\!-\!\frac{w}{\alpha}\bigg ){\mathrm d}w. 
\end{alignat*}
Let
\begin{equation*}
G^{(n,1)}_{\alpha,\beta}(s,t):=\int\limits_{\sqrt{b_{1,1}}}^{\sqrt{b_{[ns],1}}}  
\widetilde\Phi \bigg( \frac w{\alpha}\bigg ){\mathrm d}w, \qquad 
G^{(n,2)}_{\alpha,\beta}(s,t):=\int\limits_{\sqrt{b_{1,1}}}^{\sqrt{b_{1,[nt]}}}  
\widetilde\Phi \bigg( \frac w{\beta}\bigg ){\mathrm d}w.
\end{equation*}
Short calculation shows that
\begin{equation}
  \label{eq:eq3.9}
\frac 1{n^{1/2}}\Big|F^{(n,1)}_{\alpha,\beta}(s,t)-G^{(n,1)}_{\alpha,\beta}(s,t)
\Big| \leq \frac {2(\alpha\!+\!\beta)(1\!-\!\beta)}{
  \alpha\sqrt{\pi n}}\!\!\!\int\limits_{\sqrt{b_{1,1}}}^{\sqrt{b_{[ns],1}}}\!\!\!
\frac 1w {\mathrm d}w \leq \frac 1{\alpha\sqrt{\pi n}}\ln
\bigg(\frac {[ns]\!+\!1}{b_{1,1}}\bigg) \to 0
\end{equation}
as \ $n\to\infty$. \ Further, for \ $a<b$ \ we have
\begin{equation*}
\int\limits_a^b\widetilde\Phi\bigg( \frac w{\alpha}\bigg ){\mathrm
  d}w= \frac {\alpha}{\sqrt{\pi}}\left( \exp \bigg(-\frac
{b^2}{\alpha^2}\bigg)-\exp \bigg(-\frac
{a^2}{\alpha^2}\bigg)\right)+b\widetilde\Phi\bigg( \frac b{\alpha}\bigg
)- a\widetilde\Phi\bigg( \frac a{\alpha}\bigg),  
\end{equation*}
so
\begin{align*}
G^{(n,1)}_{\alpha,\beta}(s,t)=&
b_{[ns],1}^{1/2}\widetilde\Phi \bigg(\frac
  {b_{[ns],1}^{1/2}}{\alpha}\bigg )- b_{1,1}^{1/2}\widetilde\Phi
  \bigg(\frac {b_{1,1}^{1/2}}{\alpha}\bigg ) 
+\frac {\alpha}{\sqrt{\pi}}\left(\exp \bigg(-\frac
{b_{[ns],1}}{\alpha^2}\bigg)-\exp \bigg(-\frac
{b_{1,1}}{\alpha^2}\bigg)\right)
\end{align*}
that together with \eqref{eq:eq3.9} implies
\begin{equation}
   \label{eq:eq3.10}
\lim_{n\to\infty}\frac
1{n^{1/2}}F^{(n,1)}_{\alpha,\beta}(s,t)=\lim_{n\to\infty}\frac
1{n^{1/2}}G^{(n,1)}_{\alpha,\beta}(s,t)= \big (\alpha(1-\alpha)s\big)^{1/2}. 
\end{equation}
Similarly,
\begin{equation}
   \label{eq:eq3.11}
\lim_{n\to\infty}\frac
1{n^{1/2}}F^{(n,2)}_{\alpha,\beta}(s,t)=\lim_{n\to\infty}\frac
1{n^{1/2}}G^{(n,2)}_{\alpha,\beta}(s,t)= \big (\beta(1-\beta)t\big)^{1/2}. 
\end{equation}

To determine the limit of \ $n^{-1/2}F^{(n,3)}_{\alpha,\beta}(s,t)$ \
assume first that \ $(1-\beta )t<(1-\alpha )s$ \ implying \ $(1-\beta
)[nt]<(1-\alpha )[ns]$ \ if \ $n$ \ is large enough. On the one hand we have
\begin{equation*}
\frac 1{n^{1/2}}F^{(n,3)}_{\alpha,\beta}(s,t) \geq  \widetilde\Phi \bigg(\frac
  {(1-\alpha )[ns]-(1-\beta )[nt]}{\sqrt{b_{[ns],[nt]}}}\bigg ) \frac
  {\sqrt{b_{[ns],[nt]}}-\sqrt{b_{[ns],1}}}{\sqrt{n}} \to
  \sqrt{b_{s,t}}-\sqrt{b_{s,0}} 
\end{equation*}  
as \ $n\to\infty$. \ On the other hand
\begin{align*}
\frac 1{n^{1/2}}F^{(n,3)}_{\alpha,\beta}(s,t) \leq \frac
  {\sqrt{b_{[ns],[nt]}}-\sqrt{b_{[ns],1}}}{\sqrt{n}} \to
  \sqrt{b_{s,t}}-\sqrt{b_{s,0}} 
\end{align*}  
as \ $n\to\infty$, \ so 
\begin{equation}
   \label{eq:eq3.12}
\lim_{n\to\infty}\frac 1{n^{1/2}}F^{(n,3)}_{\alpha,\beta}(s,t)
=\sqrt{b_{s,t}}-\sqrt{b_{s,0}}=\big(\alpha 
    (1-\alpha)s+\beta (1-\beta)t\big)^{1/2}-\big(\alpha(1-\alpha)s\big)^{1/2}.
\end{equation}
If \ $(1-\beta )t\geq (1-\alpha )s$ \ we split the domain of
integration in \ $F^{(n,3)}_{\alpha,\beta}(s,t)$ \ into two parts,
that is \ $F^{(n,3)}_{\alpha,\beta}(s,t)=F^{(n,3,1)}_{\alpha,\beta}(s,t)
+F^{(n,3,2)}_{\alpha,\beta}(s,t)$ \ where
\begin{align*}
F^{(n,3,1)}_{\alpha,\beta}(s,t)&:=
\int\limits_{\sqrt{b_{[ns],1}}}^{\sqrt{(\alpha+\beta)(1-\beta)[ns]}}   
\widetilde\Phi \bigg(
\frac{(\alpha\!+\!\beta)(1\!-\!\alpha)[ns]}{w\beta }
-\frac{w}{\beta}\bigg ){\mathrm d}w, \\ 
F^{(n,3,2)}_{\alpha,\beta}(s,t)&:=
\int\limits_{\sqrt{(\alpha+\beta)(1-\beta)[ns]}}^{\sqrt{b_{[ns],[nt]}}}    
\widetilde\Phi \bigg(
\frac{(\alpha\!+\!\beta)(1\!-\!\alpha)[ns]}{w\beta }
-\frac{w}{\beta}\bigg ){\mathrm d}w.
\end{align*} 
Again, on the one hand we have 
\begin{align}
\frac 1{n^{1/2}}&F^{(n,3,1)}_{\alpha,\beta}(s,t)\geq
\frac 1{n^{1/2}}\!\!\!\!\!
\int\limits_{\sqrt{b_{[ns],1}}}^{\sqrt{(\alpha+\beta)(1-\beta)[ns]}}\!\!\!\!\!
\widetilde\Phi \bigg(
\frac{\sqrt{(\alpha\!+\!\beta)(1\!-\!\alpha)[ns]}-w}{\beta } \bigg
){\mathrm d}w  \label{eq:eq3.13} \\ 
=&\frac{\sqrt{(\alpha+\beta)(1-\beta)[ns]}-\sqrt{b_{[ns],1}}}{\sqrt{n}} \,
\widetilde\Phi \bigg(
\frac{\sqrt{(\alpha\!+\!\beta)(1\!-\!\alpha)[ns]}-\sqrt{b_{[ns],1}}}{\beta
} \bigg ) \nonumber \\
&+\frac {\beta}{\sqrt{n\pi}}\left(\exp \bigg(-\frac
    {\big(\sqrt{(\alpha\!+\!\beta)(1\!-\!\beta)[ns]}\!-\!
      \sqrt{b_{[ns],1}}\big)^2}{\beta  
      ^2}\bigg) \!-\!1\right)\to
  \sqrt{(\alpha\!+\!\beta)(1\!-\!\beta)s}\!-\!\sqrt{b_{s,0}}  \nonumber
\end{align}
as \ $n\to\infty$. \ On the other hand
\begin{equation*}
\frac 1{n^{1/2}}F^{(n,3,1)}_{\alpha,\beta}(s,t) \leq \frac
  {\sqrt{(\alpha\!+\!\beta)(1\!-\!\beta)[ns]}-\sqrt{b_{[ns],1}}}{\sqrt{n}} \to
  \sqrt{(\alpha\!+\!\beta)(1\!-\!\beta)s}-\sqrt{b_{s,0}} 
\end{equation*}  
as \ $n\to\infty$, \ so 
\begin{equation}
  \label{eq:eq3.14}
\lim_{n\to\infty}\frac 1{n^{1/2}}F^{(n,3,1)}_{\alpha,\beta}(s,t)
=\big((\alpha+\beta)(1-\alpha)s\big)^{1/2}-\big(\alpha(1-\alpha)s\big)^{1/2}.
\end{equation}
Similarly to
\eqref{eq:eq3.13} one can also show
\begin{align*}
-\frac 1{n^{1/2}}F^{(n,3,2)}_{\alpha,\beta}(s,t)\geq \frac 1{n^{1/2}}
\int\limits_{\sqrt{(\alpha+\beta)(1-\beta)[ns]}}^{\sqrt{b_{[ns],[nt]}}}
\widetilde\Phi \bigg(
\frac{w-\sqrt{(\alpha\!+\!\beta)(1\!-\!\alpha)[ns]}}{\beta}& \bigg
){\mathrm d}w \\
\to \sqrt{b_{s,t}}&-\sqrt{(\alpha\!+\!\beta)(1\!-\!\beta)s},
\end{align*}
and we also have
\begin{equation*}
-\frac 1{n^{1/2}}F^{(n,3,2)}_{\alpha,\beta}(s,t) \leq \frac
  {\sqrt{b_{[ns],[nt]}}-\sqrt{(\alpha\!+\!\beta)(1\!-\!\beta)[ns]}}{\sqrt{n}} \to
  \sqrt{b_{s,t}}-\sqrt{(\alpha\!+\!\beta)(1\!-\!\beta)s} 
\end{equation*}  
as \ $n\to\infty$ \ implying
\begin{equation}
   \label{eq:eq3.15}
\lim_{n\to\infty}\frac 1{n^{1/2}}F^{(n,3,2)}_{\alpha,\beta}(s,t)
=\big((\alpha+\beta)(1-\alpha)s\big)^{1/2}
-\big(\alpha(1-\alpha)s+\beta (1-\beta)t\big)^{1/2}.
\end{equation} 
Thus, by summing the limits in \eqref{eq:eq3.14} and \eqref{eq:eq3.15}
we obtain that for \ $(1-\beta )t\geq (1-\alpha )s$
\begin{equation}
   \label{eq:eq3.16}
\lim_{n\to\infty}\frac 1{n^{1/2}}F^{(n,3)}_{\alpha,\beta}(s,t)= 
2\big((\alpha\!+\!\beta)(1\!-\!\alpha)s\big)^{1/2}\!
-\!\big(\alpha(1\!-\!\alpha)s \!
  +\!\beta(1\!-\!\beta)t\big)^{1/2} \!-\!\big(\alpha(1\!-\!\alpha)s\big)^{1/2}.
\end{equation}

Finally, the asymptotic behaviour of \
$n^{-1/2}F^{(n,4)}_{\alpha,\beta}(s,t)$ \ in some sense a complementary
of the behaviour of \
$n^{-1/2}F^{(n,3)}_{\alpha,\beta}(s,t)$. \ In the same way as
\eqref{eq:eq3.12} is proved one can show that if \ $(1-\beta )t>
(1-\alpha )s$ 
\begin{equation}
   \label{eq:eq3.17}
\lim_{n\to\infty}\frac 1{n^{1/2}}F^{(n,4)}_{\alpha,\beta}(s,t)
=\big(\alpha 
    (1-\alpha)s+\beta (1-\beta)t\big)^{1/2}-\big(\beta(1-\beta)t\big)^{1/2}.
\end{equation}
In case \ $(1-\beta )t\leq (1-\alpha )s$ \ the domain of
integration in \ $F^{(n,4)}_{\alpha,\beta}(s,t)$ \ has to be split at
\ $\sqrt{(\alpha+\beta)(1-\beta)[nt]}$ \ to obtain 
\begin{equation}
   \label{eq:eq3.18}
\lim_{n\to\infty}\frac 1{n^{1/2}}F^{(n,4)}_{\alpha,\beta}(s,t)= 
2\big((\alpha\!+\!\beta)(1\!-\!\beta)t\big)^{1/2}\!
-\!\big(\alpha(1\!-\!\alpha)s\! 
  +\!\beta(1\!-\!\beta)t\big)^{1/2} \!-\!\big(\beta(1\!-\!\beta)t\big)^{1/2}.
\end{equation}
Hence, in case \ $\alpha(1-\alpha)s\leq \beta(1-\beta )t$ \ equation
\eqref{eq:eq3.8} and limits \eqref{eq:eq3.10} -- \eqref{eq:eq3.12} and
\eqref{eq:eq3.16} -- \eqref{eq:eq3.18} imply
\begin{equation}
  \label{eq:eq3.19}
\lim_{n\to\infty}\frac 1{n^{1/2}}\Var
\big(Y^{(n)}(s,t)\big)=\lim_{n\to\infty}\frac 1{n^{1/2}}
E^{(n)}_{\alpha,\beta}(s,t)=\frac {\big((1-\alpha)s\big)^{1/2} \land
  \big((1-\beta)t\big)^{1/2}}{\pi^{1/2}(\alpha+\beta)^{1/2}(1-\alpha)(1-\beta)}. 
\end{equation}

If \  $\alpha(1-\alpha)s>\beta(1-\beta )t$ \ we have
\begin{equation*}
 E^{(n)}_{\alpha,\beta}(s,t)= \frac
 1{(\alpha+\beta)(1-\alpha)(1-\beta)}
\big(E^{(n,1)}_{\alpha,\beta}(s,t)+ E^{(n,2)}_{\alpha,\beta}(s,t)+
E^{(n,3)}_{\alpha,\beta}(s,t)\big), 
\end{equation*}
with
\begin{align*}
E^{(n,1)}_{\alpha,\beta}(s,t)&:=\frac
1{2\pi}\int\limits_{b_{1,1}}^{b_{1,[nt]}}
\int\limits_{-u/\beta+(\alpha+\beta)(1-\alpha)/\beta}^{u/\alpha
  -(\alpha+\beta)(1-\beta)/\alpha} 
\frac 1u \exp\bigg(-\frac {v^2}u\bigg){\mathrm d}v{\mathrm d}u, \\
E^{(n,2)}_{\alpha,\beta}(s,t)&:=\frac
1{2\pi}\int\limits_{b_{1,[nt]}}^{b_{[ns],1}}
\int\limits_{u/\alpha-(\alpha+\beta)(1-\beta)[nt]/\alpha}^{u/\alpha-(\alpha+\beta)
  (1-\beta)/\alpha}
\frac 1u \exp\bigg(-\frac {v^2}u\bigg){\mathrm d}v{\mathrm d}u, \\
E^{(n,3)}_{\alpha,\beta}(s,t)&:=\frac
1{2\pi}\int\limits_{b_{[ns],1}}^{b_{[ns],[nt]}}
\int\limits_{u/\alpha-(\alpha+\beta)(1-\beta)[nt]/\alpha}^{-u/\beta+(\alpha+\beta)
  (1-\alpha)[ns]/\beta}
\frac 1u \exp\bigg(-\frac {v^2}u\bigg){\mathrm d}v{\mathrm d}u 
\end{align*}
and \eqref{eq:eq3.19} can be proved similarly to the other case.

Now, if \ $\alpha\beta\gamma\leq 0$ \ and $|\alpha|=1$, \
$|\beta|=|\gamma| <1$ \ or \ $|\beta|=1$, \  
$|\alpha|=|\gamma| <1$ \ using 
\eqref{eq:eq2.5} we have
\begin{equation*}
\frac 1n\Var \big(Y^{(n)}(s,t)\big)=\frac {[ns]}n \frac
{1\!-\!\gamma^{2[nt]}}{1\!-\!\gamma^2}\to \frac s{1\!-\!\gamma^2} 
\quad \text{or} \quad \frac 1n\Var \big(Y^{(n)}(s,t)\big)=\frac {[nt]}n \frac
{1\!-\!\gamma^{2[ns]}}{1\!-\!\gamma^2}\to \frac t{1\!-\!\gamma^2}, 
\end{equation*}
respectively, as \ $n\to\infty$.

At the end, if \ $\alpha=\beta=-\gamma=1$ \ the statement directly
follows from Theorem \ref{covbound}. \proofend

\markboth{\rm Baran}{\rm On the variances of a spatial
  unit root model} 

\section{Remarks on missing cases}
  \label{sec:sec4} 
\markboth{\rm Baran}{\rm On the variances of a spatial
  unit root model} 

The results of Theorem \ref{main} do not cover the cases when \
$|\alpha|<1, \ |\beta|<1, \ |\gamma|\leq 1$, \ either
$\alpha\beta\gamma <0$ \ or \ $\alpha=\beta=0$ \ is satisfied, and \
$|\alpha|-|\beta|+|\gamma|=1$ \ or \ 
$-|\alpha|+|\beta|+|\gamma|=1$ \ holds. For \ $|\gamma|<1$ \ the above
conditions yield two subcases of Case A, while for \ $|\gamma|=1$ \ we
have a subcase of Case B.

In the trivial case \ $\alpha=\beta=0$ \ and \ $|\gamma|=1$ \ using
directly \eqref{model} it is easy to see that \ $\Var
(X_{k,\ell})=k\land \ell$, \ hence
\begin{equation}
   \label{eq:eq4.1}
\lim_{n\to\infty} \frac 1n \Var \big(Y^{(n)}(s,t)\big)=s\land t.
\end{equation}

If \ $\alpha\beta\gamma <0$ \ according to the results of the Introduction it
suffices to consider \
$0<\alpha, \beta <1$ \ and \ $-1\leq \gamma <0$ \ and assume \
$\alpha-\beta-\gamma=1$ \ or \ $-\alpha+\beta-\gamma=1$. \ As the first row of
\eqref{eq:eq2.3}  holds for all positive \ $\alpha$ \ and \ $\beta$, 
\begin{equation*}
G(m,n;\,\alpha,\beta,\gamma) = \sum
                     _{r=0}^{m\land n}\binom nr \binom mr\alpha ^{m-r}\beta 
                   ^{n-r}(\alpha\beta +\gamma)^r ,                  
\end{equation*}
where
\begin{equation*}
\alpha\beta +\gamma=\begin{cases}
                     (1+\beta)(\alpha -1)<0, &\text{if \
                       $\alpha-\beta-\gamma=1$, }\\ 
                     (1+\alpha)(\beta -1)<0, &\text{if \
                       $-\alpha+\beta-\gamma=1$.} 
                    \end{cases}  
\end{equation*}
Hence, using notations of Lemma \ref{binoform}, for \
$\alpha-\beta-\gamma=1$ \ we have
\begin{equation}
  \label{eq:eq4.2}
G(m,n;\,\alpha,\beta,\gamma) = (1+2\beta)^n\sum_{r=0}^{m\land n}(-1)^r
{\mathsf P}\big (\xi_m^{(\alpha)}=m-r\big){\mathsf P}\bigg
(\eta_n^{\left( \frac{1+\beta}{1+2\beta}\right)}=r\bigg),
\end{equation}
while for \ $-\alpha+\beta-\gamma=1$ 
\begin{equation}
    \label{eq:eq4.3}
G(m,n;\,\alpha,\beta,\gamma) = (1+2\alpha)^m\sum_{r=0}^{m\land n}(-1)^r
{\mathsf P}\big (\xi_n^{(\beta)}=n-r\big){\mathsf P}\bigg
(\eta_m^{\left( \frac{1+\alpha}{1+2\alpha}\right)}=r\bigg)
\end{equation}
holds. This means that results similar to Corollary \ref{cor:cor1} can not
be obtained. Moreover, the exponential terms before the sums in
\eqref{eq:eq4.2} and \eqref{eq:eq4.3}
do not allow us to use Lemma \ref{binoappr} for
separate approximations of the probabilities behind the sums.

Finally, in case \ $\alpha=\beta<1, \ \gamma=-1$ \ short calculation
shows \citep{szego}
\begin{equation*}
G(m,n;\,\alpha,\alpha,-1)=\alpha ^{|m-n|}P_{m\land
  n}^{(0,|m-n|)}(2\alpha ^2-1),
\end{equation*}
so using notation \ $\cos (\theta)=2\alpha ^2-1$ \ we obtain
\begin{equation}
   \label{eq:eq4.4}
\Var \big (Y^{(n)}(s,t) \big)=\sum _{k=0}^{[ns]-1}\sum _{\ell=0}^{[nt]-1} \big (\cos
(\theta/2)\big )^{2|k-\ell|}\Big (P_{k\land
  \ell}^{(0,|k-\ell|)}\big(\cos(\theta)\big )\Big)^2,
\end{equation}
where \ $P_n^{(a,b)}(x)$ \ is the $n$th Jacobi polynomial with parameters
\ $a$ \ and \ $b$. \ Obviously, as \ $P_n^{(0,0)}(-1)=(-1)^n$, \ in
the trivial case \ $\alpha=\beta=0$ ($\theta=\pi$) \ limit
\eqref{eq:eq4.1} can be obtained from \eqref{eq:eq4.4}, too. However,
as in general the second parameter of the Jacobi polynomial in
\eqref{eq:eq4.4} equals \
$|k-\ell|$, \ to find the limit of the appropriately normed variances
of \ $Y^{(n)}(s,t)$ \ the
classical approximations of the Jacobi polynomials as 
e.g. Theorem 8.21.8 of \citet{szego} can not be used.

\noindent
{\bf Acknowledgments.} \ Research has been supported by 
the Hungarian  Scientific Research Fund under Grant No. OTKA
T079128/2009 and partially supported by T\'AMOP
4.2.1./B-09/1/KONV-2010-0007/IK/IT project. 
The project is implemented through the New Hungary Development Plan
co-financed by the European Social Fund, and the European Regional
Development Fund.

\end{document}